\documentclass[12pt]{article}

\usepackage{amssymb}

\usepackage{epsfig}

\makeatletter

\@addtoreset{equation}{section}

\makeatother

\newtheorem{thm}[equation]{Theorem}
\newtheorem{lemma}[equation]{Lemma}
\newtheorem{corol}[equation]{Corollary}
\newtheorem{propos}[equation]{Proposition}

\newtheorem{rmk}[equation]{Remark}
\newenvironment{rem}{\begin{rmk}\rm }{\end{rmk}}

\renewcommand\thefigure{\thesection.\@arabic\c@figure}

\renewcommand\thetable{\thesection.\@arabic\c@table}

\def\reff#1{(\ref{#1})}

\begin{document}

\def\Eq{}
\def\equ{}
\def\bn{{\bf N}}
\def\bm{{\bf M}}

\def\uV{{\underline V}}
\def\uU{{\underline U}}
\def\uu{{\underline u}}
\def\ux{{\underline x}}
\def\uw{{w_{-\infty}^{-1}}}
\def\uz{{z_{-\infty}^{-1}}}
\def\uv{{\underline v}}
\def\E{{\mathbb E}}
\def\P{{\mathbb P}}
\def\R{{\mathbb R}}
\def\Z{{\mathbb Z}}
\def\V{{\mathbb V}}
\def\N{{\mathbb N}}
\def\BB{{\mathbb B}}
\def\NN{{\bf N}}
\def\X{{\cal X}}
\def\Y{{\bf Y}}
\def\T{{\cal T}}
\def\C{{\bf C}}
\def\D{{\bf D}}
\def\G{{\bf G}}
\def\U{{\bf U}}
\def\K{{\bf K}}
\def\H{{\bf H}}
\def\J{{\cal J}}
\def\n{{\bf n}}
\def\dd{{\bf d}}
\def\b{{\bf b}}
\def\aa{{\bf a}}
\def\ee{{\rm e}}
\def\g{{\bf g}}
\def\mm{{m}}
\def\ti{{\rm TI}}
\def\sqr{\vcenter{
         \hrule height.1mm
         \hbox{\vrule width.1mm height2.2mm\kern2.18mm\vrule width.1mm}
         \hrule height.1mm}}                  
\def\square{\ifmmode\sqr\else{$\sqr$}\fi}
\def\one{{\bf 1}\hskip-.5mm}
\def\liml{\lim_{L\to\infty}}
\def\given{\ \vert \ }
\def\Given{\ \Big\vert \ }
\def\ze{{\zeta}}
\def\be{{\beta}}
\def\de{{\delta}}
\def\la{{\lambda}}
\def\ga{{\gamma}}
\def\th{{\theta}}
\def\al{{\alpha}}
\def\vep{{\varepsilon}}
\def\proof{\noindent{\bf Proof. }}
\def\A{{\bf A}}
\def\B{{\bf B}}
\def\D{{\bf D}}
\def\H{{\bf H}}
\def\h{{\bf h}}
\def\bx{{\bf X}}
\def\bz{{\bf Z}}
\def\bk{{\bf K}}
\def\bF{{\bf F}}
\def\cw{{\cal W}}
\def\zero{{\rm 0}}
\def\sign{{\rm sign}}

\def\tends#1{\mathop{\longrightarrow}\limits_{#1}}
\bibliographystyle{alpha}

\title{Processes with Long Memory: \\ Regenerative Construction and
  Perfect Simulation}
\date{}

\author{Francis Comets\footnote{Partially supported by CNRS, UMR 7599 
``Probabilit{\'e}s et Mod{\`e}les 
Al{\'e}atoires''} \\
[-3mm] 
{\it \normalsize Universit{\'e} de Paris 7}\\
Roberto Fern{\'a}ndez\footnote{
Partially supported by FAPESP, CNPq and FINEP (N{\'u}cleo de
Excel{\^e}ncia ``Fen{\^o}menos cr{\'\i} ticos em probabilidade e processos
estoc{\'a}sticos'' PRONEX-177/96)}\\
[-3mm] 
{\it \normalsize Universit{\'e} de Rouen}\\ 
Pablo A. Ferrari$^{\dagger}$\\
[-3mm] 
{\it \normalsize Universidade de S{\~a}o Paulo} \\
}  

\maketitle

\abstract{We present a perfect simulation algorithm for stationary processes
  indexed by $\Z$, with summable memory decay.  Depending on the decay, we
  construct the process on finite or semi-infinite intervals, explicitly from
  an i.i.d.\ uniform sequence. Even though the process has infinite memory,
  its value at time 0 depends only on a finite, but random, number of these
  uniform variables.  The algorithm is based on a recent regenerative
  construction of these measures by Ferrari, Maass, Mart{\'\i}nez and Ney. As
  applications, we discuss the perfect simulation of binary autoregressions
  and Markov chains on the unit interval.
}

\noindent{\bf Short Title.} Long Memory and Perfect Simulation



\endabstract
\vskip-10mm
\baselineskip 17pt 


\section{Introduction}

In this paper we consider processes with transition probabilities that depend
on the whole past history, i.e. processes with long memory.  When this
dependence decays fast enough with time, we exhibit a regenerative
construction which, besides yielding an explicit proof of existence and
uniqueness of the process, can be transcribed into a perfect simulation
scheme.  \medskip

Processes with long memory have a long history. They were first studied by
Onicescu and Mihoc (1935) under the label \emph{chains with complete
  connections} (\emph{cha{\^\i}nes {\`a} liaisons compl{\`e}tes}). Harris (1955)
proposed the somehow less used name of \emph{chains of infinite order}.
Doeblin and Fortet (1937) proved the first results on speed of convergence
towards the invariant measure. Harris (1955) called these processes chains of
infinite order, and extended results on existence and uniqueness. The chains
appeared also as part of the formalism introduced by Keane (1971 and 1976) to
study of subshifts of finite type (or covering transformations).  In this
theory the transition probabilities are called $g$-\emph{functions} and the
invariant measures $g$-\emph{measures}.

The theory of long-memory processes has found applications in the study of urn
schemes (Onicescu and Mihoc, 1935b), continued-fraction expansions (Doeblin,
1940; Iosifescu 1978 and references therein), learning processes (see, for
instance, Iosifescu and Theodorescu, 1969; Norman, 1972 and 1974), models of
gene population (Norman, 1975), image coding (Barnsley, Demko, Elton and
Gerinomo, 1988), automata theory (Herkenrath and Theodorescu, 1978), partially
observed ---or ``grouped''--- random chains (Harris, 1955; Blackwell, 1957;
Kaijser, 1975; Pruscha and Theodorescu, 1977; Elton and Piccioni, 1992) and
products of random matrices (Kaijser, 1981).  For further references we refer
the reader to Kaijser (1981 and 1994) [from which most of the material of the
previous review paragraphs is taken] and to Iosifescu and Grigorescu (1990).

It is clear that these applications should benefit from the construction and
perfect-simulation scheme presented here.  As an illustration, we discuss in
Sections \ref{ex1} and \ref{ex} applications to binary autoregressions and to
the Markov processes on the interval $[0,1]$ defined by Harris (1955) by
mapping chains with complete connection into $D$-ary expansions.

\medskip

In this paper we rely on a regenerative construction of the chain, which
generalizes, in some direction, those existing in the literature. This type of
construction has been first introduced by Doeblin (1938) for Markov chains
with countable alphabet.  Schemes for more general state spaces came much
later (Athreya and Ney, 1978; Nummelin, 1978). The first regenerative
structures for chains with complete connections were proposed by Lalley (1986,
2000) and Berbee (1987) for chains with summable continuity rates. An explicit
regenerative construction was put forward by Ferrari, Maass, Mart{\'\i}nez and
Ney (2000) in the spirit of Berbee's approach.

In the present paper we take up the scheme of Ferrari, Maass, Mart{\'\i}nez and
Ney (2000), extend it
to part of the Harris uniqueness regime and transcribe it as a perfect
simulation algorithm.  Basically, the construction used here can be
interpreted as a simultaneous coupling of \emph{all} histories, built in such
a way that at each instant $i$ there is a (random) number $k_i\ge 0$ such that
the distribution of the move $i+1$ is the same for all histories agreeing the
$k_i$ preceding instants.  This independence from the $k_i$-remote past yields
the times $\tau$ such that preceding histories are irrelevant for future
moves.  These are the regeneration times, defined by the conditions $k_i\le
i-\tau$ for all $i\ge \tau$. Both Berbee's and Lalley's constructions rely on
the regeneration probability being positive, a fact that seems to hold only
for summable continuity rates.  In contrast, our construction extends to cases
where (global) regeneration may have probability zero. Non-summable ---but
still not too slowly decreasing--- rates inside the Harris uniqueness regime
yield \emph{local} regeneration times, that is, regenerations for finite time
intervals (windows). In ergodic-theory terms, our construction is in fact a
finitary coding of a process of i.i.d.\ uniform variables in the interval
$[0,1]$.

Perfect simulation became popular after Propp and Wilson (1996) introduced the
\emph{coupling from the past} algorithm to simulate invariant measures of
Markov chains. Wilson's page {\tt
  http://dimacs.rutgers.edu/\char126dbwilson/exact} provides updated and
extensive references on perfect simulation. Foss and Tweedie (1998) and
Corcoran and Tweedie (1999) proposed a general framework based on regeneration
schemes for Markov chains and the so-called ``stochastic recursive
sequences''. These are processes defined by $X_{n+1} = f(X_n,\xi_n)$, where
$\xi$ is a stationary process. In the Markovian case $\xi_i$ are i.i.d.; in the
non-Markovian case, it remains the matter of how to construct/simulate the
sequence $\xi_i$. Our algorithm applies to a wide variety of non-Markovian
processes and, through the formalism of random systems with complete
connections, it can be used to simulate Markov processes with state space of
large cardinality (e.g.\ the unit interval). Section \ref{ex}
present an example along these lines.

Our main result is an explicit construction, as a deterministic function of a
sequence of i.i.d. random variables uniformly distributed in $[0,1]$, of
realizations of the stationary chain with infinite memory.  As corollaries we
get: (i) an alternative proof of the existence and uniqueness of the
stationary process, and (ii) a perfect simulation algorithm and a regeneration
scheme for this process. These results are summarized in Theorem \ref{a9} and
Corollaries \ref{z1}, \ref{50} and \ref{a58}. In Section \ref{def1} we
introduce the basic definitions and in Section \ref{ex1} we 
illustrate the simulation on a concrete example,
which is continued in
Section \ref{ex}. The results of Section \ref{s.res} are proved in
Sections 
\ref{s.fin}-\ref{a57}
and \ref{s.regener}, though the perfect  simulation algorithm is described in
Section \ref{a56}  in its general version.

\section{Definitions}\label{def1}
We denote by $G$ our alphabet, $\N^*=\N\setminus\{0\}$ and $-\N^*
=\{-i:i\in\N^*\}$.  In what follows $G$ can be finite or countable, though in
the latter case conditions \reff{n.1} and \reff{15} below impose severe
limitations. The set $G^{-\N^*}$ is the space of \emph{histories}; we write
$\omega_a^b :=(w_b, w_{b-1},\ldots,w_a)$ for $-\infty\le a \le b \le +\infty$.
For shortness we write $\underline w = \omega_{-\infty}^\infty$. Let
$P:G\times G^{-\N^*}\to[0,1]$ be a probability transition kernel; that is,
$P(g\vert \uw)\ge 0$ for all $g$ and
\begin{equation}
  \label{a1}
  \sum_{g\in G} P(g\vert \uw) = 1\,.
\end{equation}
for each $\uw \in G^{-\N^*}$. The kernel $P$ defines, by telescopic
products, what in statistical mechanics is called a \emph{specification}. A
specification is a consistent system of conditional probabilities, where
consistency is required for \emph{all} histories $\uw$.  In standard
probabilistic treatments, such requirements are made only almost surely with
respect to some pre-established appropriate measure.  But in the present
setting the determination of the appropriate measure is part of the problem,
and stronger requirements are a priori necessary.

Denoting by $\eta(i), i \in \Z,$ the coordinate mappings on $G^\Z$, we
say that a (non-necessarily stationary) probability measure $\nu$ on
$G^\Z$ ---or a \emph{process} with distribution $\nu$--- 
is \emph{compatible} with the specification $P$ 
if the latter
is a version of the one-sided conditional probabilities of the former:
 \begin{equation}
   \label{a4}
   \nu\Bigl(\eta\,:\, \eta(i)=g \;\Bigm|\; 
\eta(i+j) = w_j,\, j\in -\N^*\Bigr) =
   P(g\vert \uw) 
 \end{equation}
for all $i\in \Z$, $g\in G$ and $\nu$-a.e. $\uw\in
G^{-N^*}$. Then the identities
\begin{equation}
\label{f5}
\nu\Bigl(\eta\;:\; \eta(i+k)=g_k, k \in [0,n]\; \Bigm|\; \eta(i+j)=w_j,\,
  j\in -\N^*  \Bigr)
 = \prod_{k=0}^n P(g_k\vert g_{k-1}, \dots g_0, \uw)\,,  
\end{equation}
hold for $\nu$-a.e. $\uw$, where the concatenation is defined
by 
$$
(w_{k},w_{k-1},\dots, w_{\ell}, {\uz}) :=
(w_{k},w_{k-1},\dots, w_{\ell}, z_{-1}, z_{-2},\dots)\;.
$$
The regenerative construction of this paper is based on the
functions
\begin{eqnarray}
a_0(g) &:=& \inf\Bigl\{P(g\vert {\uz}) \ : \  
{\uz} \in G^{-\N^*} \Bigr\}\nonumber\\[5pt]
a_k(g \vert w_{-k}^{-1}) &:=&
\inf\Bigl\{P(g\vert w_{-1},\dots, w_{-k}, {\uz}) \ : \ 
{\uz} \in G^{-\N^*} \Bigr\}\;, \quad k\ge 1\,,
\label{n.2}
\end{eqnarray}
defined for $k \in\N$, $g, w_{-1},\ldots, w_{-k}\in G$.  [These
functions are denoted $g(i_0|i_{-1},\ldots,i_{-k})$ by Berbee (1987)].
The numbers
\begin{equation}
  \label{a2}
 a_k := \inf_{w_{-k}^{-1}} \ 
\sum_{g\in G} a_k(g \vert w_{-k}^{-1}) \;\;\;,
\end{equation}
$k \in \N$, determine a probabilistic threshold for memories limited
to $k$ preceding instants.  The sequences $a_k(g \vert w_{-k}^{-1})$
and $a_k$ are non-decreasing in $k$ and contained in $[0,1]$.

For our construction we shall use a sequence $\uU=(U_i\,:\,i\in\Z)$
of independent random variables with uniform distribution in $[0,1[$,
constructed on the corresponding canonical probability space
$(\Omega,{\mathcal F},\P)$. We denote $\E$ the expectation with respect to
$\P$. 

\section{An example: binary autoregressive processes}\label{ex1}

To motivate the method, we present an example that shows how to construct and
perfect simulate a process with infinite memory. Let us consider binary
autoregressive processes. Such a process is the binary version of
autoregressive (long memory) processes used in statistics and econometrics.
It describes binary responses when covariates are historical values of the
process (McCullagh and Nelder, 1989, Sect. 4.3).

Let the state space be $G=\{-1,+1\}$, $\th_0$ a real number 
and $(\th_k; k \geq 1)$ a summable real sequence. Let $q: \R \mapsto ]0,1[$
be strictly increasing and continuously differentiable. Assume that 
\begin{equation}
  \label{ex01}
P(\, \cdot\,  \vert \uw) {\rm \ is \ the\  Bernoulli\  law\
  on\  }\{-1,+1\} {\rm \ with \ parameter } \;\;
q\Bigl(\th_0 +\sum_{k \geq 1} \th_k w_{-k}\Bigr)\;\;,
\end{equation}
i.e., $P( + 1 \vert \uw)= q(\th_0 +\sum_{k \geq 1} \th_k w_{-k})=
1-P( - 1 \vert \uw)$. 

By compactness there exists at least one process compatible with \reff{ex01};
the conditions for uniqueness are well known. The question is how to simulate
(construct) such a stationary process. To do that we construct a family of
partitions of the interval $[0,1[$ indexed by $k$ and $w_{-k}^{-1}$ using
$a_k(1| w_{-k}^{-1})$ and $a_k$ defined in the previous
section. Letting 
$$r_k=
\sum_{m>k} |\th_m|
$$ 
we have
\begin{eqnarray}
  \label{pk2}
 a_k(1| w_{-k}^{-1})&=&q \left( \th_0 + \sum_{1 \leq m \leq k} \th_m w_{-m}-
   r_k\right) \nonumber \\
a_k(-1|w_{-k}^{-1})&=&
1-q \left( \th_0 + \sum_{1 \leq m \leq k} \th_m w_{-m} +r_k \right)
\end{eqnarray}
and
\begin{equation}
  \label{ex011}
a_k\;=\; 1- \sup_{ w_{-k}^{-1}} \left\{q \left( \th_0 + \sum_{1 \leq m
      \leq k} \th_m w_{-m}  
+r_k \right)
-q \left( \th_0 + \sum_{1 \leq m \leq k} \th_m w_{-m}-
      r_k\right)\right\} \;.
\end{equation}
For each $k$ and $w_{-k}^{-1}$, let
$\B_k(g|w_{-k}^{-1})$ be intervals of length
$a_k(g|w_{-k}^{-1})-a_{k-1}(g|w_{-k+1}^{-1})$ placed consecutively in
lexicographic order in $k$ and $g$, starting at the origin. Use the convention
$a_0(g|w_{0}^{-1})= a_0(g)$ and $a_{-1} = 0$. These intervals form a partition
of the interval $[0,1[$.

Let $(U_i\,:\,i\in\Z)$ be the sequence of i.i.d. random variables uniformly
distributed in $[0,1[$ defined at the end of Section \ref{def1}.  For
$n\in\Z$ define the random variable
\begin{equation}
  \label{ea13}
  K_n := \sum_{k\ge 0} k \;\one\{U_n\in [a_{k-1},a_k)\}\,.
\end{equation}
In our construction the variable $K_n:=K_n(\uU)$ indicates how many sites in
the past are needed to compute the state at time $n$. To each $n$ associate an
arrow going from $n$ to $n-K_n$. The state at site $n$ will be independent of
the states at $t\le s$ if no arrow starting at $\{s,\dots,n\}$ finishes to the
left of $s$. Let $\tau[n]:=\tau[n](\uU)$ be the maximum of such $s$:
\begin{eqnarray}
  \label{x13}
  \tau[n]&:= &\max\Bigl\{s \leq n\,:\, U_j < a_{j-s},\,
  j\in[s,n]\Bigr\}\nonumber\\ 
&=& \max\Bigl\{s \leq n\,:\, K_j \ge s,\,
  j\in[s,n]\Bigr\} 
\end{eqnarray}
Notice that $\tau[n]$ is a stopping time for the sequence $(U_{n-k}:k\ge 0)$.
We show in Theorem~\ref{a9} that the condition
\begin{equation}
  \label{ex010}
\sum_k r_k < \infty
\end{equation}
is sufficient to guarantee $\P(\tau[n]>-\infty)=1$ which, in turn, is 
an equivalent condition to the feasibility of the following construction.

\paragraph{Simulation (construction) of the stationary measure}

\begin{enumerate}
\item Generate successively i.i.d. uniform random variables $U_n, U_{n-1},
  \dots$. Stop when $U_{\tau[n]}$ is generated. Using \reff{ea13}, the values
  $K_n, K_{n-1}, \dots, K_{\tau[n]}$ are simultaneously obtained.
\item Use the $U_{\tau[n]},\dots, U_n$ and $K_{\tau[n]},\dots, K_n$ generated
  in the previous step to define
  \begin{equation}
    \label{paf1}
  X_{j}\;=\;g \qquad\hbox{ if }\qquad
U_{j}\in \bigcup_{\ell=0}^{K_j} \B_\ell(g|X_{j-\ell}^{j-1}) \;,
  \end{equation}
 recursively from $j=\tau[n]$ to $j=n$.
\item Return $X_n$. The algorithm has also constructed
  $X_{\tau[n]},\dots,X_{n-1}$.
\end{enumerate}
The expression \reff{paf1} is well defined because by the definition of
$\tau[n]$, $j-K_j\ge \tau[n]$, 
$U_j\in [0, a_{K_j}[ \subset 
\bigcup_{g \in G}  \bigcup_{\ell=0}^{K_j} \B_\ell(g|X_{j-\ell}^{j-1})$ 
and the set in \reff{paf1}
depends at most on $X_{j-K_j},\dots, X_{j-1}$. This is discussed in detail in
Section~\ref{a57}; see~\reff{a21}. We show in Theorem \ref{a9} that the above
algorithm constructs a realization $X_n$ of a random variable which has the
one-coordinate marginal of the unique measure compatible with the
specification.  To construct a realization of the measure in a finite window,
just repeat the algorithm for other $n$ reusing always previously generated
$U_j$. The algorithm induces a function $\Phi: [0,1[^\Z \to \{-1,1\}^\Z$.
Given the event $\{\tau[n] = k\}$, $X_n$ ---the $n$th coordinate of
$\Phi(U_{-\infty}^{\infty})$--- depends only on $U_k, \dots,U_n$. (To be
rigorous, one should give the definition of $\Phi$ when $\tau[n]=-\infty$ for
some $n$; this is an arbitrary but irrelevant choice as this set has
probability zero under our hypotheses.)

This construction exploits a well known fact. The existence of a renovating
event gives rise to a perfect simulation algorithm and a regenerative
structure. In our case there is a regeneration at time $s$ if $j-K_j\ge s$ for
all $j\ge s$ (no arrow passes over $s$). However for the construction of the
measure in site $n$ we use a weaker condition: it suffices that no arrow
passes from $j$ to the left of $s$ for $j\in[s,n]$.

While in some cases one may have explicit expressions for $a_k$, in general,
this will not be the case (see \reff{ex011}). 
An useful aspect of our construction is that in
these cases we can work with {\bf lower bounds} $a_k^*$. We shall discuss this
issue in Section \ref{ex}.

\section{Results}\label{s.res}

The existence of our regeneration scheme depends on the non-increasing
sequence
\begin{equation}
  \label{b8}
 \be_m\,:=\, \prod_{k= 0}^m a_k \;,
\end{equation}
with $a_k$ defined in \reff{a2}, and a necessary (but not sufficient)
condition is $\lim_{k \to \infty} a_k=1$. 
The regeneration time for a window
$[s,t]$, for $-\infty<s<\infty$ and $s \le t \le \infty$, is the random
variable
\begin{equation}
  \label{13}
  \tau[s,t]:= \max\Bigl\{m \leq s\,:\, U_k < a_{k-m}\,,
  k\in[m,t]\Bigr\} 
\end{equation}
which may be $-\infty$.  Notice that $\tau[s,t] = \min\{\tau[n]\,:\,
n\in[s,t]\}$ and that it is a \emph{stopping time} for $(U_{t-k}\,:\,
k\ge 0)$, in the sense that $ \{\tau[s,t]\ge j\} \in \mathcal F(U_i:
i\in[j,t])$ for $j\le s$, where $\mathcal F (U_i: i\in[j,t])$ is the sigma
field generated by $(U_i: i\in[j,t])$.

Our main result is:  

\begin{thm}
  \label{a9} 
If
\begin{equation}
\sum_{m\ge 0}\be_m=\infty
\label{n.1}
\end{equation}
then 
  \begin{itemize}
  \item [(i)] For each finite $[s,t]\subset\Z$,
  $\P(\tau[s,t]>-\infty)=1$, where $\tau[s,t]$ is defined in \reff{13}.
\item [(ii)] There exists a measurable function
  $\Phi:[0,1[^\Z\to G^\Z$ ---described in Section \ref{a57}---
  such that
  \begin{equation}
    \label{a25}
\mu\; :=\; \P(\Phi(\uU)\, \in \,\cdot\,)\,,
  \end{equation}
  the law of $\Phi(\uU)$, is compatible with~$P$. Moreover, the distribution
  $\mu$ is stationary.

\item[(iii)] In addition, the function $\Phi$ has the property that
  for each finite interval $[s,t]\subset\Z$, its restriction
  $\{\Phi(\uU)(i)\,:\,i\in[s,t]\}$ depends only on the
  values of $U_i$ in the interval $[\tau[s,t],t]$. More precisely, for
  $i\in [s,t]$,
\begin{equation}
  \label{a12}
  \Phi(\uU)(i) = \Phi(\, \dots, v_{\tau
  -2},v_{\tau-1},U_\tau,\dots, U_t,v_{t+1},\dots)(i)
\end{equation}
for any sequence $\underline v\in[0,1[^\Z$ (we abbreviated $\tau[s,t]$
as $\tau$).
\item[(iv)] The law of $\tau[s,t]$ satisfies the following bound:
  \begin{equation}
    \label{55}
    \P(s-\tau[s,t]>m)\,\le\, \sum_{i=0}^{t-s} \rho_{m+i}
  \end{equation}
where $\rho_m$ is the probability of return to the origin at epoch
$m$ of the Markov chain on $\N$ starting at time zero at the origin
with transition probabilities
\begin{eqnarray}
p(x,x+1) &=& a_x \nonumber\\
p(x,0) &=& (1-a_x) 
\label{n.3}
\end{eqnarray}
and $p(x,y)=0$ otherwise.  In particular, if $(1-a_k)$ decreases
exponentially fast with $k$, then so does $\rho_k$. If $(1-a_k) \sim
k^{-\gamma}$ with $\gamma>1$ 
then $\rho_k$ decreases with the same power. 

\item[(v)] If
\begin{equation}
  \label{15}
  \be:=\lim_{m\to\infty}\be_m >0
\end{equation}
then items (i), (iii) and (iv) above hold also for $t=\infty$.
\end{itemize} 
\end{thm}

This theorem is proven in Sections \ref{s.fin} and \ref{a57}.  More detailed
bounds on the parameters $\rho_m$ are given in Proposition \ref{60}. 
Conditions \reff{n.1} and \reff{15} require $a_0>0$ as in Harris~(1955). Both
conditions are noticeable weaker than those imposed by Lalley~(1986),
Berbee~(1987) and Bressaud, Fern{\'a}ndez and Galves~(1999b), as well as those
corresponding to the $g$-measure approach (Ledrappier, 1974) and to Gibbsian
specifications (Kozlov, 1974).

\begin{rem}\label{frem.f}
  The construction can be performed replacing the
$a_k$'s with lower bounds, that is, considering a sequence $a_k^*$ in
$]0,1[$, such that
\begin{equation}
\begin{array}{ll}
\displaystyle
a_k^* \leq a_k  & \hbox{for } k \geq 0 \;.\\
\end{array}
  \label{frem0*1r}
\end{equation}
The theorem is valid replacing the unstarred $a_k$ by starred ones, and
using the corresponding starred versions of $\beta_m$, $\tau$ and $\rho_m$.
While the actual $a_k$'s give shorter regeneration times, they could be hard
to estimate. Suitable choices of $a^*_k$ could provide a reasonable compromise
between shorter regeneration times and feasible calculational efforts (Section
\ref{ex}).
\end{rem}

Our first corollary is the uniqueness of the measure compatible with
$P$.

\begin{corol}[Loss of memory and uniqueness]
  \label{z1}
\mbox{}

\begin{itemize}
\item[(i)] Every measure $\mu$ compatible with the specification $\P$
  has the following loss-of-memory property: If $f$ is a function
  depending on the interval $[s,t]$ and $i \leq s$,
\begin{equation}
\label{loss.f}
\Bigl| \mu\Bigl(f\Bigm|\eta(j)=w_j,\, j<i\Bigr) - 
\mu\Bigl(f\Bigm|\eta(j)=v_j,\, j<i\Bigr) \Bigr| 
\;\le\;  2\,\|f\|_\infty\, \sum_{j=0}^{t-s} \rho_{s+j-i}
\end{equation}
for every $\underline w, \underline v \in G^\Z$. Here we use the
  notation $\mu f := \int \mu(d\eta)f(\eta)$.
\item[(ii)] If $\sum_{m\ge 0}\be_m=\infty$ the measure $\mu$ defined
  in \reff{a25} is the unique measure compatible with~$P$. 
\end{itemize}
\end{corol}

The uniqueness result is not new.  Under the more restrictive condition
\reff{15} it was already obtained by Doeblin and Fortet (1937).  Harris (1955)
[see also Section 5.5 of Iosifescu and Grigorescu (1990) and references
therein] extended this uniqueness to a region that coincides with~\reff{n.1}
for two-symbol alphabets but it is larger for larger alphabets.  Other
uniqueness results, in smaller regions, were obtained by Ledrappier (1974) and
Berbee (1987) in different ways.  Results on loss of memory were also obtained
by Doeblin and Fortet (1937) (see also Iosifescu, 1992), under the summability
condition \reff{15}.  Bressaud, Fern{\'a}ndez and Galves (1999b) extended those
to a region defined by a condition slightly stronger than \reff{n.1}.  The
rates of loss of memory obtained in this last references strengthen those of
Doeblin, Fortet and Iosifescu, but are weaker than ours. 

A corollary of (ii) and (iii) of the theorem is a perfect simulation scheme:

\def\wPhi{{\widetilde \Phi}}
\begin{corol}[Perfect simulation]
  \label{50}
  Let $P$ be a specification with $\sum_{m\ge 0}\be_m=\infty$ and
  $\mu$ the unique measure compatible with $P$. 
For each finite window  $[s,t]$ there exist
a family   $(\wPhi_{m,t}; m \leq s)$ of functions $\wPhi_{m,t}\,
:\,[0,1[^{[m,t]} \to G^{[m,t]}$
  such that
\begin{equation}
  \label{51}
  \mu \Bigl(\eta\,:\,\eta(i)= g_i,\; i\in [s,t]\Bigr)\;=\; 
  \P\Bigl(\uU\,:\,\wPhi_{\tau,t}(U_{\tau},\dots,U_t)(i) 
  = g_i,\; i\in [s,t]\Bigr)
\end{equation}
where $\tau=\tau[s,t]$ is the stopping time defined in \reff{13}.
\end{corol}

Expression \reff{51} is our perfect simulation scheme for $\mu$.
Possible implementations are discussed in Section \ref{a56}. Roughly
speaking the algorithm goes as follows: 
\begin{enumerate}
\item Produce a realization of $U_t,\dots,U_s,U_{s-1}, \dots,
  U_{\tau[s,t]}$; 
\item Compute $\wPhi_{\tau[s,t],t}(U_{\tau[s,t]},\dots,U_t)$.
\item Return the configuration
  $(\wPhi_{\tau[s,t],t}(U_{\tau[s,t]},\dots,U_t)(i)\,:\, i\in [s,t])$.  Its
  law is the projection of $\mu$ in the window $[s,t]$.
\end{enumerate}

Our algorithm shares two features with the famous Coupling From The Past
(CFTP) algorithm for Markov chains, introduced by Propp and Wilson (1996): (1)
The r.v. $U_i$ are generated sequentially backwards in time, until the
stopping time $\tau[s,t]$ is discovered.  (2) The algorithm is based on a
coupled realization of the chain for all possible initial conditions ensuring
the coalescence of the different trajectories before the observation window.
Restricted to the Markovian case, our renewal times are, in principle, larger
than the coalescence times of well designed CFTP algorithms.  Nevertheless,
our approach, besides extending to non-Markovian processes, yields coupling
times that depend only on the numbers $a_k$ defined in \reff{a2} [see
\reff{imp.f}], and hence that are only indirectly related to the cardinality
of the alphabet.

The bounds in (iv) of Theorem \ref{a9} can be used to control the
\emph{user-impatience bias}.  This picturesque name proposed by Fill (1998)
relates to the bias caused by practical constraints such as time limitation or
computer storage capacity which force the user to stop a (long and unlucky)
run before it reaches perfect equilibrium.  Indeed, suppose that while
sampling window $[s,t]$ we decide to abort runs with regeneration times
$\tau[s,t] >M$ for some large $M>0$, causing our algorithm to produce a biased
distribution $\widehat\mu_{[s,t]}^M$.  Applying Lemma 6.1 of Fill (1998) and
\reff{55}, we obtain:
\begin{equation}
\sup_{A\in G^{[s,t]}} \, \Bigl| \mu(A) \,-\, \widehat\mu_{[s,t]}^M(A)
\Bigr| \;\le\; {\displaystyle \sum_{i=0}^{t-s} \rho_{_{\scriptstyle
      M+i}}  \over \displaystyle 1 -
\sum_{i=0}^{t-s} \rho_{_{\scriptstyle M+i}}}\;.
\label{imp.f}
\end{equation}
\medskip

In regime (v), Theorem \ref{a9} yields the following
\emph{regeneration} scheme. Let $\NN\in\{0,1\}^\Z$ be the random
counting measure defined by
\begin{equation}
  \label{a59}
  \NN(j) := \one\{\tau[j,\infty]=j\}\,.
\end{equation}
Let $(T_\ell\,:\,\ell\in\Z)$ be the ordered time events of $\NN$
defined by $\NN(i)=1$ if and only if $i=T_\ell$ for some~$\ell$,
$T_\ell<T_{\ell+1}$ and 
$T_0\le 0 < T_1$. 

\begin{corol}[Regeneration scheme]
\label{a58}  
If $\be>0$, then the process $\NN$ defined in \reff{a59} is a
stationary renewal process with renewal distribution
\begin{equation}
\P(T_{\ell+1}-T_{\ell}\geq m) \;=\; \rho_m
\label{ff.f}
\end{equation}
for $m>0$ and $\ell \neq 0$.  Furthermore, the random vectors
$\xi_\ell\in\cup_{n\ge 1} G^n$, $\ell\in\Z$, defined by
\begin{equation}
  \label{a60}
 \xi_\ell(\uU)\; 
 :=\;(\Phi(\uU)(T_\ell),\dots,\Phi(\uU)(T_{\ell+1}-1)) 
\end{equation}
are mutually
independent and $(\xi_\ell(\uU)\,:\, \ell\neq 0)$ are
identically distributed. 
\end{corol}


\section{Distribution of $\tau[s,t]$ and  $\tau[s,\infty ]$}
\label{s.fin}

In this section we prove items (i) and (iv) [and the corresponding
parts in (v)] of the theorem.  We fix an increasing sequence of
numbers $(a_k)$ such that
\begin{equation}
a_k \nearrow 1\ \hbox{as } k \nearrow \infty
\label{nea.f}
\end{equation}
and define the following  \emph{house of cards}
family of chains $((W^m_n\,:\,n\ge m) \,:\,m\in\Z)$ by $W^m_m = 0$ and
for $n>m$:
\begin{equation}
  \label{1}
  W^m_n := (W^m_{n-1} +1)\, \one\{U_n<a_{{\scriptstyle W^m_{n-1}}}\} 
\end{equation}
where $U_n$ are the uniform random variables introduced in the previous
section. \emph{Notation alert:} Here $W^m_n$ is not a vector of lenght $n-m$,
but the position at time $n$ of a chain starting at time $m$ at the origin.
Alternatively, $\P(W^m_n=y\,|\,W^m_{n-1}=x)=p(x,y)$, where the latter
are the transition probabilities defined in item (iv) of Theorem
\ref{a9}.  Notice that 
\begin{equation}
  \label{rho.f}
  \rho_k \;=\; \P(W^m_{m+k}=0)
\end{equation}
for all $m\in\Z$ and $k\in \N^*$.  The  monotonicity of the
$a_k$'s implies that
\begin{equation}
W^m_n \;\geq\; W^k_n \quad \hbox{for all } m<k \leq n\;.
\label{mon.f}
\end{equation}
Hence, $W^m_n=0$ implies that
$W^k_n=0$ for $m<k\leq n$, and the chains \emph{coalesce} at
time $n$:
\begin{equation} 
   \label{coalesce}
W^m_n=0 \; \Longrightarrow\; W^m_t=W^k_t, \; m \leq k \leq n \leq t
\end{equation}
(furthermore, $W^m_t > W^k_t$ for $t$ smaller than
the smallest such $n$).

By the definition (\ref{13}) of $\tau$, we have for $j\leq s$
\begin{eqnarray} 
  \tau[s,t] < j 
  & \iff & \forall m \in [j,s], \exists n \in [m,t] : W^{m-1}_n=0
  \nonumber\\
  & \iff & \forall m \in [j,s], \exists n \in [s,t] : W^{m-1}_n=0
  \nonumber\\
  & \iff & \max\{ m<s : \forall n \in [s,t] ,  W^{m}_n>0\}<j-1
  \nonumber\\
  & \iff & \exists n \in [s,t] : W^{j-1}_n=0
  \label{f3}
\end{eqnarray}
where the second line follows from the coalescing property
\reff{coalesce}, and the last line is a consequence of the second one
and the monotonicity \reff{mon.f}. From the third line we get
\begin{equation} 
   \label{3}
  \tau[s,t] \;=\; 1\,+\, \max \Bigl\{m< s: W^m_n >0, \,\forall\,
  n\in[s,t]\Bigr\}\;. 
\end{equation}

\begin{lemma}  \label{4}\mbox{}

\begin{itemize}
\item[a)] Condition $\sum_{n\ge 0} \prod_{k=1}^n a_k =\infty$
is  equivalent to any of the two following  properties:
\begin{itemize}
\item[(a.1)] The chain $(W^m_n\,:\, n\ge m)$ is non positive-recurrent for
each $m$.
\item[(a.2)] For each $-\infty<s\le t<\infty$, $\tau[s,t]>-\infty$ \emph{a.s.}.
\end{itemize}

\item[b)] Condition $\;\prod_{k=1}^\infty a_k \;>\; 0\;$
is  equivalent to any of the two following  properties:
\begin{itemize}
\item[(b.1)] $(W^m_n\,:\, n\ge m)$ is transient
\item[(b.2)] For each $-\infty<s$, $\tau[s,\infty]>-\infty$ \emph{a.s.}.
\end{itemize}
\end{itemize}

In both cases, inequality \reff{55} holds.

\end{lemma}

\proof a) It is well known that $W_n$ is positive-recurrent if and
only if $\sum_{n\ge 0} \prod_{k=1}^n a_k \;<\;\infty$, (for instance
Dacunha-Castelle, Duflo and Genon-Catalot (1983), p.~62 ex.~E.4.1).
Also, by \reff{f3}
\begin{equation}
  \label{6}
  \Bigl\{\tau[s,t]< m\Bigr\} = \bigcup_{i\in[s,t]}
\Bigl\{W^{m-1}_i=0\Bigr\}
\end{equation}
for $m\le s$.  By translation invariance, the probability of the
right-hand side
of \reff{6} satisfies
\begin{equation}
  \label{7}
  \P\Bigl(\bigcup_{i\in[s,t]}\{W^0_{-m+i}=0\}\Bigr)\;\in \; 
\left[ \P(W^0_{t-m} =0)\;,\, \sum_{i=1}^{t-s}
  \P(W^0_{s-m+i} =0) \right]
\end{equation}
where the right-hand-side is a consequence of the monotonicity
property \reff{mon.f}.  As $m\to -\infty$ this interval remains
bounded away from 0 in the positive-recurrent case, but shrinks to 0
otherwise.  Since $\P(\tau[s,t]=-\infty) = \lim_{m\to
  -\infty}\P(\tau[s,t]< m)$, this proves that (a.1) and (a.2) are
equivalent.  \smallskip

\noindent b) Clearly,
\begin{equation}
  \label{9}
  \P(W^m_n \neq 0\,,\;\forall n\ge m)\;=\; \prod_{i=0}^\infty a_i\;.
\end{equation}
This implies that the product of $a_k$ is positive if and only if the
house-of-cards process is transient.  From
\reff{6} we have that
\begin{equation}
  \label{10}
  \P(\tau[s,\infty]<m) = \P\Bigl(\bigcup_{i\in[s-m+1,\infty]}\{W^0_i=0\}\Bigr)
\end{equation}
which goes to zero as $m\to-\infty$ in the transient case only. Therefore
(b.1) and (b.2) are equivalent.
\smallskip

Inequality \reff{55} follows from \reff{6} or \reff{10}, due to \reff{rho.f}.
\square \medskip

The following proposition is due to Bressaud, Fern{\'a}ndez and Galves (1999b).

\begin{propos}
  \label{60} 
  Let $a_k$ be a $[0,1]$-sequence increasing to one. Let $\rho_k$ be the
  probability of return to the origin at epoch $k$ of $(W^0_n\,:\,n\ge 0)$.
\begin{itemize}
\item[(i)] If  $\sum_{n\ge 0} \prod_{k=1}^n a_k \;=\;\infty$, then
  $\rho_n\to 0$.

\item[(ii)] If $\prod_{k=1}^\infty a_k \;>\; 0$, then
$\sum_{n\ge 0} \rho_n <\infty$.

\item[(iii)] If $(1-a_n)$ decreases exponentially then so does $\rho_n$.

\item[(iv)] 
  \begin{equation}
    \label{61}
  \prod_{k=1}^\infty a_k \;>\; 0 \quad\hbox{ and }\quad \sup_i\limsup_{k\to\infty}\left({1-a_i\over1-a_{ki}}\right)^{1/k}\;\le\;1
    \quad\hbox { imply }\quad \rho_n = {\rm O}(1-a_n)
  \end{equation}

\end{itemize}
\end{propos}

Item (iv) can be applied, for instance, when $a_n\sim 1- (\log n)^b
n^{-\gamma}$ for $\gamma>1$. Items (i) and (ii) are direct transcriptions of
(a.1) and (b.2) of the previous lemma. 
\medskip

\begin{rem}\label{frem3.r}
  The results of this section are also monotonic on the choice of
  sequence $(a_k)$: A sequence $(a^*_k)$, satisfying \reff{nea.f},
  with $a^*_k\le a_k$ for all $k$, yields chains $W^{*m}_n$ with
  larger probability of visiting $0$, and hence with larger
  regeneration intervals $s-\tau^*[s,t]$ and larger values of
  $\rho^*_m$. This is the content of Remark \ref{frem.f}.
\end{rem}

\section{Construction of $\Phi$}\label{a57}

In this section we prove results needed to show (ii) and (iii) [and the
corresponding part in (v)] of Theorem \ref{a9}. The results hold for any
sequence $(a^*_k)$ satisfying \reff{frem0*1r}, but for notational simplicity
we omit the superscript ``$*$''.  For $g\in G$ let $b_{0}(g) \,:=\, a_{0}(g)$,
and for $k\ge 1$,
$$ 
b_k(g \vert w_{-k}^{-1}) \,:=\, a_k(g \vert w_{-k}^{-1})
-a_{k-1}(g \vert w_{-k}^{-1}).
$$
For each $\uw\in G^{-\N^*}$ let $\{\B_k(g \vert w_{-k}^{-1})\,:\,
g\in G,\, k\in\N\}$ be a partition of $[0,1[$ with the following
properties: (i) for $g\in G$, $k\ge 0$, $\B_k(g \vert w_{-k}^{-1})$
is an interval closed in the left extreme and open in the right one with
Lebesgue measure $|\B_k(g \vert w_{-k}^{-1})|
=b_k(g \vert w_{-k}^{-1})$; (ii) these intervals
are disposed in increasing lexicographic order with respect to $g$ and
$k$ in such a way that the left extreme of one interval coincides with
the right extreme of the precedent:
$$ 
\B_0(g_1),\B_0(g_2),\dots, 
\B_1(g_1|w_{-1}^{-1}),\B_1(g_2|w_{-1}^{-1}),
\dots,
\B_2(g_1|w_{-2}^{-1}),\B_2(g_2|w_{-2}^{-1}),\dots
$$
is a partition of [0,1[ into consecutive intervals increasingly arranged.
This definition is illustrated in Figure~1  in the case $G=\{1,2\}$.

\begin{figure}
            \centerline{\input{figfn.pstex_t}}
        \centerline{\small{\it Figure 1}: Partition 
$\{\B_k(g \vert w_{-k}^{-1})\,:\,
g\in G,\, k\in\N\}$ of $[0,1[$ in the case $G=\{1,2\}$}
\end{figure}

In particular we have
\begin{equation}
  \label{a23}
 \Bigl|\bigcup_{k\ge 0} \B_k(g \vert  w_{-k}^{-1}
)\Bigr|\,=\,\sum_{k\ge 0} |\B_k(g \vert w_{-k}^{-1})|\,=\,P(g|\uw)
\end{equation}
and
\begin{equation}
  \label{a24}
\bigcup_{g\in G}\,\bigcup_{k\ge 0} \B_k(g  \vert  w_{-k}^{-1})
\,=\,[0,1[
\end{equation}
where all the unions above are disjoint.
For $k \ge 0$ let
$$
\B_k(  w_{-k}^{-1}) :=
\bigcup_{g\in G}  \B_k(g \vert  w_{-k}^{-1}).
$$

By \reff{a2} and \reff{frem0*1r}, we have
\begin{equation}
  \label{a21}
  [0,a_k[\;\;\subset \;\;\bigcup_{\ell=0}^k
\B_\ell( w_{-\ell}^{-1}) = \left[0 \;,\, \sum_{g \in G} 
a_{k}(g \vert { w_{-k}^{-1}}) \right[ \;\;
, \qquad\hbox{ for all }\uw \in G^{-\N^*}.
\end{equation}
As a consequence, 
\begin{equation}
[0,a_k[\,\cap\, \B_\ell(w_{-\ell}^{-1}) \;=\; \emptyset \quad
\hbox{for } \ell > k\;,
\label{empty.f}
\end{equation}
a fact that makes the definitions
\begin{eqnarray}  \label{34}
\B_{\ell,k}(w_{-1},\dots,w_{-k} ) &:=& [0,a_k[ \,\cap\,
\B_\ell(w_{-\ell}^{-1}) \\
\B_{\ell,k}(g \vert w_{-1},\dots,w_{-k}) &:=& [0,a_k[ 
\,\cap\, \B_\ell(g \vert w_{-\ell}^{-1}) \;.
\label{34.1}
\end{eqnarray}
meaningful \emph{for all} $\ell, k\ge 0$.  

Items (ii) and (iii) [and the corresponding part in (v)] of Theorem
\ref{a9} are immediate consequences of the following proposition. 
Similar to \reff{x13}, define 
$$
 \tau[n](\uu):= \max\Bigl\{s \leq n\,:\, u_j < a_{j-s},\,
 j\in[s,n]\Bigr\}\;.
$$

\begin{propos}
  \label{27}
Let $g_0$ be an arbitrary point in $G$.
Define the function 
$\Phi:[0,1[^\Z\to G^\Z , \uu \mapsto \ux=\Phi(\uu)$,
recursively from $j=\tau[n](\uu)$ to $j=n$ in the following fashion:
write  $\tau=\tau[n](\uu)$ and define
\begin{eqnarray}
 x_{\tau} &:= &\sum_{g\in G} g \, \one\Bigl\{ u_{\tau} 
\in \B_{0}(g) \Bigr\} \nonumber\\
 x_{\tau+1} &:= &\sum_{g\in G} g \, \one\Bigl\{ u_{\tau+1} 
\in \B_{0}(g) \,\cup\, \B_{1,1}(g \vert x_{\tau})\Bigr\} 
\nonumber\\
& \vdots& \nonumber\\
 x_{n} &:= &\sum_{g\in G} g \, \one\Bigl\{ u_{n} \in \bigcup_{0 \leq j
   \leq \tau}
\B_{j,n-\tau}(g \vert x_{n-1},\dots,x_{\tau}) \Bigr\} \; 
  \label{a27}
\end{eqnarray}
if $\tau[n](\uu)$ is finite, and $x_n=g_0$ if $\tau[n](\uu)=-\infty $. (Each
sum in \reff{a27} reduces to a single term.)  Then:
\begin{itemize}
\item[(i)] $\Phi$ is well defined and measurable.
\item[(ii)] For each $n\in\Z$, the component $\Phi(\uu)(n)$ depends
  only on the values of $u_i$ on the interval $[\tau[n](\uu),n]$ if
  $\tau[n](\uu)>-\infty$, and on $(u_i:i\le n)$ otherwise.
\item[(iii)] If $\P(\tau[n](\uU)>-\infty)=1$, for all
  $n\in\Z$, then the law $\mu$ of $\Phi(\uU)$ is compatible
  with $P$ and $\mu$ is stationary.
\end{itemize}
\end{propos}

\proof 
Let
\begin{equation}
A_n \;=\;\Bigl\{\uu \in [0,1]^\Z :\tau[n](\underline
u)>-\infty\Bigr\}\;.
\label{an.f}
\end{equation}


\noindent {\sl (i)}  On the set $A_n$, the consistency of the definition
\ref{a27} follows from two facts: (1) By definition,
$\tau[n]\le \tau[j]$ for $j\in[\tau[n],n]$, and (2) \reff{a21} shows
that if the event $\{U_n < a_k\}$ holds, then we only need to look at
$x_{n-1},\dots,x_{n-k}$ to obtain the value of $x_n$.  These facts
imply that for every $k\in[\tau[n],n]$, the value of $x_k$ computed
using \reff{a27} with $k$ in place of $n$ yields the same value as the
one obtained as part of the recursive calculation \reff{a27} for
$x_n$.  The $\mathcal F(U_i: i\leq n)$-measurability of $x_n$ follows,
then, from definition \reff{a27}.  As the sets $A_n$ are $\mathcal
F(U_i: i\leq n)$-measurable, see \reff{13}, the maps $x_n$ remain
measurable after gluing. We conclude that $\Phi(\uu)=(x_n; n
\in \Z)$ is well defined and measurable.
\smallskip

\noindent {\sl (ii)}  Immediate from the definition (now known to be
consistent). 
\smallskip

\noindent {\sl (iii)}  
When $A_n$ is true, \reff{empty.f} implies that Definition \ref{a27}
amounts to
\begin{equation}
 x_{n} \,= \,\sum_{g\in G} g \, \one\Bigl\{ u_{n} \in \bigcup_{k \geq 0 }
\B_{k}(g \vert x_{n-k}^{n-1}) \Bigr\} \;.
\label{xn.f}
\end{equation}
which implies
\begin{eqnarray}
\P \Bigl(\Phi(\uU)(n)=g \Bigm\vert U_i: i< n\,;\, A_n\Bigr)
\;=\; \P \Bigl(U_{n} \in \bigcup_{k \geq 0 }
\B_{k}(g \vert [\Phi(\uU)]_{n-k}^{n-1})\Bigm\vert U_i: i< n\,;\,
A_n\Bigr)
\label{xn2.f}
\end{eqnarray}
Since each $\Phi(\uU)(j)$ is $\mathcal
F(U_i: i\leq j)$-measurable, and each $U_n$ is independent of $\mathcal
F(U_i: i< n)$, 
if $A_n$ has full measure, then \reff{xn2.f} equals 
\begin{equation}
  \label{ab23}
 \left|\bigcup_{k\ge 0} \B_k\Bigl(g \Bigm\vert [\Phi(\uU)]_{n-k}^{n-1}
\Bigr)\right|\,=\,\sum_{k\ge 0} \Bigl|\B_k\Bigl(g \Bigm\vert [\Phi(\uU)]_{n-k}^{n-1}\Bigr)\Bigr|\,=\,P\Bigl(g\Bigm|[\Phi(\uU)]_{-\infty}^{n-1}\Bigr)
\end{equation}
by \reff{a23}. In other words,
\begin{equation}
\P \Bigl(\Phi(\uU)(n)=g \Bigm\vert U_i: i< n\Bigr) \;=\;
 P\Bigl(g \Bigm\vert \Phi(\uU)(i): i <n\Bigr) \;.
\end{equation}
As the right-hand side depends only on $(\Phi(\uU)(i), i <n)$, so does the
left-hand side. Hence the law of $\Phi(\uU)(n)$ given $(\Phi(\uU)(i), i <n)$
itself is still given by $P$.  Therefore the distribution $\mu$ of $\Phi(\uU)$
is compatible with $P$ in the sense of \reff{a4}. It is stationary by
construction. \square \medskip


\begin{rem}
  \label{frem0*} 
  The previous argument shows, in particular, that if $a_k^* \leq a_k$
  for all $k \geq 0$, and $\tau^*[s,t] (\uu)$ is finite (in
  which case $\tau[s,t]<\tau^*[s,t]$ is finite, see Remark
  \ref{frem3.r}), then $\Phi (\uu) (j)= \Phi^* (\uu)
  (j)$ for all $j \in [s,t]$.  (Each $\Phi^*(\uu)(j)$ depends
  on a larger number $j-\tau^*[j]$ of preceding $U's$.)  If Condition
  \reff{n.1} holds for the coefficients $\beta_m^*= \prod_{k=0}^m
  a_k^*$, both $\Phi^*$ and $\Phi$ have laws compatible with $P$.
\end{rem}
\noindent

\vskip 3mm
\noindent{\bf Proof of Corollary \ref{z1}.}
We represent the one-sided conditioned measure through a family of
functions $\Phi(\cdot|\underline w,i):[0,1]^\Z\mapsto G^\Z$
corresponding to a history $\underline w\in G^\Z$ frozen for times
smaller than $i\in\Z$.  Writing, for shortness, $y_j=\Phi(\underline
U|\underline w,i)(j)$, we set $y_j=w_j$ for $j< i$ and successively
for $n= i, i+1, \ldots$,
\begin{equation}
  \label{b27}
 y_{n} \,:= \,\sum_{g\in G} g \, \one\Bigl\{ U_{n} \in \bigcup_{k\ge 0}
\B_k \Bigl(g \Bigm\vert y_{n-1},\dots,y_i,w_{i-1},\dots,
w_{n-k},\dots\Bigr) \Bigr\}\;.  
\end{equation}

Let $\mu$ be any probability measure on $G^\Z$ compatible with $P$.
{From} \reff{f5} and \reff{a23} we see that the law of $\Phi(\underline
U|\underline w,i)$ is a regular version of $\mu$ given
$\eta_{i-1}=\omega_{i-1}, \eta_{i-2}=\omega_{i-2},\dots$.  That is,
\begin{equation}
\mu\Bigl(f\Bigm |\eta(j)=w_j,\, j<i\Bigr) \;=\;
\E[f(\Phi(\uU|\uw,i))] 
\label{cond.f}
\end{equation}
for any continuous $f$.  We now follow the classical arguments.
\smallskip

\noindent{\sl Proof of (i):}  Let $f$ be as stated. By \reff{cond.f}
\begin{equation}
\mu\Bigl(f\Bigm|\eta(j)=w_j,\, j<i\Bigr) - 
\mu\Bigl(f\Bigm|\eta(j)=v_j,\, j<i\Bigr) 
\;=\;  \E\Bigl[f(\Phi(\uU|\uw,i)) - f(\Phi(\uU|v_{- \infty}^{-1},i))\Bigr] \;.
\label{loss2.f}
\end{equation}
Since
\begin{equation}
  \label{22}
  \one\Bigl\{\Phi(\uU|\underline w,i)(n)\ne \Phi(\uU|\underline
  v ,i)(n)\Bigr\}\;\le\; \one\{\tau[n]\le i\} 
\end{equation}
the absolute value of both terms in \reff{30} is bounded above by
\begin{equation}
  \label{31}
  2\,\|f\|_\infty\, \P(\tau[s,t]<i)\;.
\end{equation}
To conclude, we use the bound \reff{55}.
\smallskip

\noindent{\sl Proof of (ii):} If $\mu$ and $\mu'$ are two measures on
$G^\Z$ compatible with $P$, 
\begin{eqnarray}
  \label{30}
  |\mu f - \mu' f| &=& \Bigl|\int \mu(d \uw) \,
\mu\Bigl(f\Bigm |\eta(j)=w_j,\, j<i\Bigr)\,-\, 
\int \mu'(d \uv)\, \mu'\Bigl(f\Bigm|\eta(j)=v_j,\, j<i\Bigr)\Bigr|
  \nonumber\\
&\le& \int\int \mu(d \uw)\, \mu'(d \uv)\,
\Bigl| \mu\Bigl(f\Bigm|\eta(j)=w_j,\, j<i\Bigr) - 
\mu\Bigl(f\Bigm|\eta(j)=v_j,\, j<i\Bigr) \Bigr| \;,\nonumber\\
\ 
\end{eqnarray}
which, by part (i), goes to zero as $i\to-\infty$. \square

\paragraph{Proof of Theorem \protect\ref{a9}}

We finish this section showing how Theorem \protect\ref{a9} follows
from previous results.

Lemma \ref{4} proves (i) and (iv) and the corresponding part of (v). The
bounds mentioned at the end of (iv) are consequence of Proposition \ref{60}.
Proposition \ref{27} proves parts (ii) and (iii) and the corresponding part in
(v).

\section{Perfect simulation}\label{a56}

In this section we prove Corollary \ref{50}. The construction of the function
$\widetilde\Phi_{\tau,t}$ relies on an alternative construction of the
stopping time $\tau[s,t]$. 
Assume $s\le t<\infty$ and define
\begin{equation}
  \label{a44}
  Z[s,t]:= \max\{K_n-n+s\,:\, n\in [s,t]\} \geq 0
\end{equation}
where $K_n = K_n(\uU)$ is defined in \reff{ea13}. $Z[s,t]$ is the number of
sites to the left of $s$ we need to know to be able to construct $\Phi$ in the
interval $[s,t]$. Let $Y_{-1} := t+1$, $Y_0:= s$ and for $n\ge 1$,
inductively
\begin{eqnarray}
  \label{a14}
  Y_n &:=& Y_{n-1}- Z[Y_{n-1},Y_{n-2}-1]
\end{eqnarray}
Then it is easy to see that
\begin{equation}
  \label{a15}
  \tau[s,t] \,=\, \lim_{n\to\infty} Y_n\, =\,\max \{Y_n\,:\, Y_n=Y_{n+1}\}
\;\; {\rm a.s.},
\end{equation}
with the usual convention $\max \emptyset = -\infty$.

\paragraph{ Construction of the perfect-simulation function
  $\widetilde\Phi_{\tau,t}$} 

\begin{enumerate}
\item Set $Y_{-1}= t+1$, $Y_0=s$ and iterate the following step
until $Y_n=Y_{n-1}$:
\begin{itemize}
\item 
Generate $U_{Y_{\ell}},\dots,U_{Y_{\ell-1}-1}$. Use \reff{ea13} to compute
  $K_{Y_{\ell}},\dots,K_{Y_{\ell-1}}$  and \reff{a44} and \reff{a14}
  to compute $Y_{\ell+1}$.
\end{itemize}
\item Let $\tau = Y_n$
\item Iterate the following procedure from $k=\tau$ to $k=t$:
  \begin{itemize}
  \item \label{step5} Define $x_k$ using \reff{a27}
  \end{itemize}
\item Return $(x_j\,:\, j\in
  [s,t])\;=\;\widetilde\Phi_{\tau,t}(U_\tau,\dots,U_t)(j)\,:\, j\in [s,t])$.
\end{enumerate}
By definition of $\tau$, $\widetilde\Phi_{\tau,t}(U_\tau,\dots,U_t)(j)=
\Phi(\uU)(j)$, $j\in [s,t]$. \square

\begin{rem}
  \label{frem*} 
  The algorithm can be applied for any choice of $(a_k^*, k \geq 0)$
  satisfying \reff{frem0*1r} and, in addition,
$$
\sum_m \prod_{k=0}^m a_k^* = \infty \;\;.
$$
The smaller the $a^*_k$, the smaller the stopping times
$\tau^*[s,t]$ of the resulting perfect-simulation scheme.  Also the
return probabilities $\rho^*_m$ increase if the $a^*_k$ increase,
worsening the bound \reff{imp.f} on the user-impatience bias.
\end{rem}

\section{Regeneration scheme}\label{s.regener} 
In this section we prove Corollary \ref{a58}. The stationarity of
$\NN$ follows immediately from the construction.  Let
\begin{equation}
  \label{36}
  f(j)\, := \,\P\Bigl(\NN (-j)=1\,|\, \NN(0)=1\Bigr)
\end{equation}
for $j\in\N^*$.  To see that $\NN$ is a renewal process it is sufficient
to show that 
\begin{equation}
  \label{35}
  \P\Bigl(\NN(s_\ell) = 1\,;\, \ell=1,\dots,n\Bigr)\;=\; 
\be \prod_{\ell=1}^{n-1} f(s_{\ell+1}-s_\ell) 
\end{equation}
for arbitrary integers $s_1<\cdots<s_k$. 
[From Poincar{\'e}'s inclusion-exclusion formula, a measure on
$\{0,1\}^\Z$ is characterized by its value on cylinder sets of the
form $\{\zeta\in\{0,1\}^\Z\,:\, \zeta(s)=1,\, s\in S\}$ for all finite
$S\subset \Z$.  For $S= \{s_1,\dots,s_k\}$, a renewal process must
satisfy \reff{35}.]  For $j\in\Z$, $j'\in\Z\cup\{\infty\}$, define
\begin{eqnarray}
  \label{38}
  H[j,j'] \,:=\,\cases{  \{U_{j+\ell}< a_{\ell},
\ell=0,\dots,j'-j\},&if $j\le j'$\cr
&\cr
\hbox{``full event''},&if $j>j'$ }
\end{eqnarray}
With this notation,
\begin{equation}
  \label{39}
\bn(j) = \one \{H[j,\infty]\},\ \ j\in\Z. \label{500}
\end{equation}
and
\begin{equation}
  \label{42}
  \P\Bigl(\NN(s_\ell) = 1\,;\, \ell=1,\dots,n\Bigr)
\;=\;\P\Bigl\{\displaystyle \bigcap_{\ell=1}^n H[s_\ell,\infty]\Bigr\}
\end{equation}
{From} monotonicity we have  for $j<j'<j''\le \infty$,
\begin{equation}
  \label{40}
  H[j,j'']\cap H[j',j'']\,\, =\,\, H[j,j'-1]\cap H[j',j''],
\end{equation}
and then, with $s_{n+1}=\infty$ we see that \reff{42} equals
\begin{eqnarray}
  \label{41}
 \,\prod_{i=1}^n\P\Bigl\{H[s_\ell,s_{\ell+1}-1]\Bigr\}\;,
\label{42b}
\end{eqnarray}
On the other hand,
\begin{equation}
  \label{ng1}
  f(j) = \P(H[-j,\infty]\,|\, H[0,\infty])\,=\,\P(H[-j,-1]) 
\end{equation}
Hence, \reff{42b} equals the right hand side of \reff{35} and we have proved
that $\NN$ is a renewal process.  \smallskip

On the other hand, by stationarity,
\begin{equation}
\P (T_{\ell +1}-T_{\ell} \geq m) 
\;=\; \P\Bigr(\tau[-1, \infty]<-m+1 \Bigm\vert \tau[0, \infty]=0\Bigr)
\end{equation}
and, hence, by \reff{f3} and \reff{rho.f}
\begin{equation}
\P (T_{\ell +1}-T_{\ell} \geq m) 
\;=\; \P( W_{-1}^{-m+1}=0) \;=\; \rho_m \;\;,
\end{equation}
proving \reff{ff.f}.
\smallskip

The independence of the random vectors
$\xi_\ell$ follows from the definition of $T_\ell$ and part (iii) of
Theorem \ref{a9}.  \square

\section{Applications}\label{ex}

\subsection{Binary autoregressions, continued}\label{ex.1}

In this subsection we continue the discussion of example \reff{ex01}. 
Recalling the notations of Section \ref{ex1}, we define 
\begin{eqnarray}
C^+ &=& \max\Bigl\{ q'(x): x \in \bigl[\th_0-\sum_{m>0} |\th_m|
\,,\,  \th_0+\sum_{m>0} |\th_m|\bigr]\Bigr\}\\
C^- &=& \min\Bigl\{ q'(x): x \in \bigl[\th_0-\sum_{m>0} |\th_m|
\,,\, \th_0+\sum_{m>0} |\th_m|\bigr]\Bigr\}\;.
\end{eqnarray}
From
Definition \reff{15}, a simple computation shows that
\begin{equation}
  \label{fex010}
\be >0 \;\Longleftrightarrow\; \sum_k r_k < \infty 
\;\Longleftrightarrow \;\sum_k k\, |\th_k| < \infty
\end{equation}
and also that for $ |\th_k| \sim Ck^{-2}$, Condition \reff{n.1} is
satisfied for $C<(2C^+)^{-1}$, but not satisfied for $C>(2C^-)^{-1}$.
Hence $\limsup_k k^2 |\th_k| <(2C^+)^{-1} $ is sufficient
for being in the Harris regime~\reff{n.1}, and $\limsup_k
|\th_k|k^2 \ln^2k < \infty $ implies that $\be >0$.

Since $a_k$ given by \reff{ex011} has no simple expression in general,
Remark~\ref{frem.f} could be useful. Indeed,
\begin{equation}
  \label{fex02}
1-2C^+r_k \;\leq\; a_k \;\leq\; 1-2C^- r_k\;.
\end{equation}

Under the extra condition that $\th_k \neq 0$ for infinitely many
$k$'s, we have in fact $ a_k \sim 1-2C^+ r_k $ as $k \to \infty$. 

We can replace the coefficients $a_k$ with the following lower bounds.
We choose first some integer $k_0$ such that $2C^+r_{k_0}<1$ and define
\begin{equation}
\begin{array}{ll}
\displaystyle
a_k^*=a_k \wedge (1-2C^+r_{k_0}) & {\rm for} \; k < k_0\;,\\
\displaystyle
a_k^*=1-2C^+r_{k} & {\rm for} \; k \geq k_0\;.
\end{array}
\label{extra.f}
\end{equation}
We can use the modification of our algorithm at the end of Section \ref{a56}
with these coefficients $a_k^*$. Note that we only need to compute at most
$k_0$ different $a_k$'s.
 
We now discuss two well-studied choices for $q$.

\paragraph{Sigmo{\"\i}d case:} In addition we assume here that $q$ is concave
on $\R^+$ with $q(x)+q(-x)=1, x \in \R$. One natural choice is
\begin{equation}
  \label{fex03}
q(x)= \frac{  \exp x}{2 \cosh x}= \frac{  1}{2(1+\exp(-2x))}\;\;,
\end{equation}
i.e., the so-called logistic function and logit model (Guyon (1995), Ex.
2.2.4), where the explicative variables are the values of the process in the
past.  For a general sigmo{\"\i}d $q$, the supremum in \reff{ex011} is
achieved for $ w_{-k}^{-1}$ minimizing $|\th_0 + \sum_{1 \leq m \leq k} \th_m
w_{-m}|$.

\paragraph{Linear case:} We take now $q(x)=(1+x)/2$, and necessarily 
$|\th_0| +\sum_{m>0} |\th_m|<1$. As we will see, linearity makes the
model \reff{ex01} somehow trivial, but also instructive.
Writing ${\cal B}(p)$ for the 
Bernoulli distribution with parameter $p$ and $\delta$ the Dirac measure,
we note that the one-sided conditional law \reff{ex01} is given by
a convex combination
\begin{eqnarray}
P( \cdot  \vert \uw) &=& {\cal B}\left(
\frac{1 + \th_0+\sum_{k\geq 1} \th_k w_{-k}}{2}\right)= {\cal B}\left(
\frac{1+ \th_0-r_0}{2} + \sum_{ k\geq 1} |\th_k|\frac{\sign (\th_k)
  w_{-k}+1}{2}\right)\nonumber\\ 
 &=& (1-r_0)\,{\cal B}\Bigl((1+ \th_0-r_0)/2(1-r_0)\Bigr)+
\sum_{ k\geq 1} |\th_k|\, \delta_{\sign (\th_k) w_{-k}}\;\;,
\label{ld.f}
\end{eqnarray}
since  ${\cal B}(\lambda p+ (1- \lambda) p')=\lambda {\cal B}(p)+
 (1- \lambda) {\cal B}(p')$ for $\lambda \in [0,1]$. 
In this example we have $a_k(w_{-k}^{-1})=1-r_k=a_k$
independent of $w$, and for $k \geq 1$, 
\begin{equation}
b_k(\pm 1\vert w_{-k}^{-1})\;=\;(\th_k w_{-k}+ |\th_k|)/2 
\;=\; \left\{\begin{array}{ll}
|\theta_k| & \hbox{if }w_{-k}= \sign(\theta_k)\\
0 & \hbox{otherwise}
\end{array}
\right.
\ \ \in \{0, |\th_k|\}\;.
\end{equation}
Hence one of the two intervals $\B_k(\pm 1 \vert w_{-k}^{-1})$ is empty, while
 the other one $\B_k(\sign (\th_k)w_{-k} \vert w_{-k}^{-1})$ has
 length $ |\th_k|$.
This is in accordance with formula \reff{ld.f}.
 
In other respects, the decomposition \reff{ld.f} can be directly interpreted
in terms of simulation: the value 
of the process at time $i=0$ is chosen according to a ``new'' coin tossing with 
probability $1-r_0$, and set to the value
$\sign (\th_k) w_{-k}$ with probability $|\th_k|\; (k=1,2, \ldots)$.
In the latter case the value of $ w_{-k}$ is needed, and will be constructed using
\reff{ld.f} again, etc$\ldots$ Clearly this recursive construction
will eventually stop if and only if $|\th_0| +\sum_{m>0} |\th_m|<1$. However, 
our construction in this paper requires
the extra condition \reff{n.1}, which loosely speaking, amounts
to $\limsup_k k^2 |\th_k|<1$. The reason for stronger assumptions is that, 
in order to cover
general processes, we need in our general construction to check {\it all
intermediate times} between 0 and $\tau{[0]}$, though in the special case 
of $P$ given by \reff{ld.f}, it is not necessary to construct all of them
following the above lines.


\subsection{Markov chains defined by $D$-ary expansions}\label{9.2}

These are processes having the unit interval as ``alphabet'', $I=[0,1]$,
and defined through another, auxiliary, process with a finite
alphabet.  Formally, a family of maps is established between sequences
of a finite alphabet $G=\{0,1,\cdots,D-1\}$ and real numbers in $I$
via $D$-ary expansions: For each $n\in\Z$
\begin{equation}
  \label{har1.f}
  \begin{array}{rcl}
X_n: G^\Z & \longrightarrow & I\\
(\eta(i): i\in\Z) & \mapsto & x_n=\sum_{j=1}^\infty
\eta(n-j)/D^j \;.
\end{array}
\end{equation}
This map induces a natural map from probability kernels $P:G\times
G^{-\N^*}\mapsto[0,1]$ to probability kernels $F:I\times
I \mapsto[0,1]$: For each $x\in I$, given an $\uw \in
G^{-\N^*}$ with $x=X_0(\uw)$
\begin{equation}
  \label{har2.f}
  F\Bigl(X_1={g+x\over D} \Bigm| X_0=x\Bigr) \;=\; P(g|\underline
  w)\;.
\end{equation}
Interest focuses on the existence and properties of measures on the
Borelians of $I^\Z$ compatible with such a 
(Markov) 
kernel $F$.

Maps \reff{har1.f}--\reff{har2.f} have been already introduced by
Borel (1909) for i.i.d.\ $\eta(i)$.  The general case in which the
$\eta(i)$ form a chain with long memory is the object of Harris (1955)
seminal paper.  They are the prototype of the random systems with
complete connections mentioned in the Introduction.  Harris determines
conditions for the existence and uniqueness of these processes,
through the study of long-memory chains: if the finite-alphabet chain
satisfies a condition similar to (but weaker than) \reff{n.1}, there
is a unique process $\mathcal M$ on $I^\Z$ compatible with $F$.
This process $\mathcal M$ is of course a (stationary) Markov chain with 
transition probability kernel $F$. Harris shows that 
its marginal distribution is continuous, except in the degenerate case
with constant $\eta(i)$'s where it is
concentrated on one of the points $0,
1/(D-1), 2/(D-1), \cdots, 1$.
Furthermore, if the process is mixing and not degenerate, 
this marginal is purely singular whenever the variables 
$\eta(i)$ are not
independent uniformly distributed (in which case the marginal is uniform).


Our approach yields, in a straightforward way, a perfect-simulation
scheme for the measures $\mathcal M$ obtained in this fashion, if the
auxiliary process $\eta(i)$ satisfies condition \reff{n.1}.  Indeed,
the map 
\begin{equation}
  \label{har4.f}
  \begin{array}{rcl}
X: G^{-\N^*}  &\longrightarrow& I\\
\uw & \mapsto & x=X_0(\uw)
\end{array}
\end{equation}
can be made bijective by fixing rules to decide between sequences which are
eventually 0 and those that are eventually $D-1$.  In turns, this map induces
a bijection between the sigma algebra $\mathcal S_\ell$ formed by unions of
intervals with endpoints in multiples of $D^{-\ell}$, and the subsets of
$G^{[-\ell,-1]}$.  We conclude that, if $\mu$ and $\mathcal M$ are the
processes compatible with the kernels $P$ and $F$ related as in \reff{har2.f},
then the restriction of $\mathcal M$ to $\mathcal S_\ell$ can be perfectly
simulated by mapping, via $X$, the perfect samples of the measure $\mu$ on the
window $[-\ell,-1]$, obtained by the algorithm of Section \ref{a56}.  We point
out that the union of the families of $\mathcal M_\ell$ uniquely determines
the measure $\mathcal M$ (it forms a so-called $\pi$ system).

\paragraph{Conclusion} 
For Markov chains in general state-space with transition kernel
$P(x,\,\cdot\,)$ satisfying the Doob's condition $P(x,\,\cdot\,)\ge
\beta\varphi(\,\cdot\,)$, for all state $x$, some $\beta>0$ and a measure
$\varphi$ on the state-space, the forward coupling is well known. The
corresponding coupling-from-the-past algorithm is illustrated in Example 2 of
Foss and Tweedy (1999) and in Corcoran and Tweedie (1999). Notice however that
the mere existence of a minorization measure is not sufficient to construct
the couplings: one needs to explicitely know $\varphi$ and $\beta$.  In this
section we have discussed an example of a Markov chain with state-space
$[0,1]$ that can be transcribed as a chain with complete connections and
state-space $\{0,1, \ldots D-1\}$. Perfect simulating the latter provides a perfect
simulation for the former. Instead of the exibition of a minorization measure
for the Markov chain our method requires the 
knowledge of $a_k(g|\uw)$ and lower bounds of
$a_k$ for the related infinite-memory chain.

\section*{Acknowledgements}
RF wants to thank Davide Gabrielli and Jeff Steif for useful discussions.
FC acknowledge Xavier Guyon and Alexander Tsybakov for pointing out
references relating to Section \ref{ex.1}. 

FC and RF thanks IME-USP for hospitality. PAF thanks the Laboratoire des
Probabilit{\'e}s de l'Universit{\'e} de Paris VII, the Departement de
Mathematiques de l'Universit{\'e} de Cergy Pontoise and the Laboratory Rapha{\"e}l
Salem of the Universit{\'e} de Rouen, for hospitality during the completion of
this work.

We thank one of the anonymous referees
for his careful reading of the manuscript and for his useful suggestions that
improved the paper.

Francis Comets:
Universit{\'e} Paris 7 -- Denis Diderot, Math{\'e}matiques, case 7012, 2
place Jussieu, 
75251 Paris Cedex 05, France.
Email: comets@math.jussieu.fr\\

Roberto Fernandez: Universit{\'e} de Rouen,           
CNRS-UPRES-A 6085                
Math{\'e}matiques, site Colbert     
F 76821 Mont Saint Aignan - France
Email: Roberto.Fernandez@univ-rouen.fr\\

Pablo Ferrari: Universidade de S{\~a}o Paulo, Instituto de Matem{\'a}tica
e Estat{\'\i}stica, Caixa Postal 66281, 05315-970 S{\~a}o Paulo, Brazil.
Email: pablo@ime.usp.br

\end{document}